\documentclass[preprint]{article}
\usepackage{amsmath}
\usepackage{amsthm}
\usepackage{amsfonts}
\usepackage[margin=1in]{geometry}
\usepackage{amssymb}
\usepackage{array}
\usepackage{tikz}
\usepackage{float}
\usepackage{xcolor}
\usepackage{graphicx}
\usepackage{framed}
\usepackage{todonotes}
\usepackage{algorithm}
\usepackage{algpseudocode}
\usepackage{listings}
\usepackage{color}
\usepackage{comment}
\usepackage{nicefrac}

\definecolor{mygreen}{RGB}{28,172,0} 
\definecolor{mylilas}{RGB}{170,55,241}

\newtheorem{theorem}{Theorem}[section]
\newtheorem{remark}{Remark}[section]
\newtheorem{corollary}{Corollary}[theorem]
\newtheorem{lemma}[theorem]{Lemma}
\newtheorem{definition}[theorem]{Definition}

\renewcommand{\em}{\it}
\usepackage{xcolor}  

\usepackage{subfig}
\usepackage{graphicx}
\usepackage{amsfonts,amsmath,yfonts}
\usepackage{amssymb}

\usepackage{subcaption}  

\usepackage{subcaption} 

\newcommand{\field}[1]{\mathbb{#1}}
\newcommand{\Rept}{{\rm Re}}

  \newcommand{\bi}{\begin{itemize}}
  \newcommand{\ei}{\end{itemize}}
  \newcommand{\be}{\begin{enumerate}}
  \newcommand{\ee}{\end{enumerate}}
  \newcommand{\beq}{\begin{equation}}
  \newcommand{\eeq}{\end{equation}}
  \newcommand{\beqa}{\begin{eqnarray}}
  \newcommand{\eeqa}{\end{eqnarray}}

\numberwithin{equation}{section}

\title{Convergence Analysis of Optimal SOR for a Class of Consistently Ordered 2-Cyclic Matrices with Complex Spectra\thanks{The work of the second author was funded in part by Discovery grant number RGPIN-2023-05244 from the
Natural Sciences and Engineering Research Council of Canada (NSERC).}}

\date{July 23, 2025}

\author{L. Robert Hocking\thanks{Department of Computer Science, The University of British Columbia,  \texttt{\{rhocking,greif\}@cs.ubc.ca}}
\and Chen Greif\footnotemark[2]}

\begin{document}

\maketitle

\begin{abstract}
\vspace{0.15cm}
\noindent
Asymptotic rates of convergence of optimal SOR applied to linear systems with consistently ordered 2-cyclic matrices have been extensively studied in the case where the Jacobi eigenvalues are are real and contained in an interval centered at the origin.  It is well known that as the rightmost endpoint of the interval approaches $1$ from below,  optimal SOR converges an order of magnitude faster than Jacobi.  
We generalize this to the situation where the Jacobi spectrum is contained in a line segment in the complex plane that is symmetric about the origin.  This is an important class of linear systems, which arise often in various physical applications; complex-shifted linear systems are included in this family. Optimal relaxation parameters are known in this case, but a detailed convergence analysis  does not seem to exist in the literature. Using techniques of complex analysis, we derive convergence rates, finding that in the complex case they are affected not only by the distance to 1
of the right-hand endpoint of the line segment 
as in the real case, but also by its phase.
\end{abstract}

\vskip 2mm
\noindent
{\bf Keywords.} iterative methods for linear systems, successive overrelaxation (SOR), optimal relaxation parameter, convergence analysis, complex analysis, complex shift  

\vskip 2mm
 \noindent 
 {\bf Mathematics Subject Classification.} 65F10, 65N22

\section{Introduction}

The successive overrelaxation (SOR) method  is a classical stationary iterative method for solving large and sparse linear systems of the form
\begin{equation} 
A x= b.
\label{eq:system}
\end{equation}
We assume that the matrix $A$ is $n \times n$ and can be complex, and the vectors $x$ and $b$ are of length $n$.
Let us define the splitting
\begin{equation}
A=D-E-F,
\label{eq:splitting}
\end{equation} 
where $D$ is the diagonal of $A$, assumed nonsingular, and $E$ and $F$ the negations of its strictly lower triangular and upper triangular parts, respectively. 
Given an initial guess, $x_0$, and a scalar parameter $\omega$ (referred to henceforth as ``the relaxation parameter'') the SOR iteration is given by
\begin{subequations}\label{eq:SOR}      
\begin{align}
  x_{k+1} &= \mathcal{L}_\omega x_k + (D-\omega E)^{-1}\,\omega b,
            \qquad k = 0,1,\dots                           \label{eq:SORiter}\\
  \intertext{where}                                   
  \mathcal{L}_\omega &= (D-\omega E)^{-1}\!
                        \bigl(\omega F + (1-\omega)D\bigr).
                                                       \label{eq:SORitermat}
\end{align}
\end{subequations}
The matrix $\mathcal{L}_\omega$ is known as the {\em SOR iteration matrix}.

A beautiful convergence analysis for SOR exists when $A$ is a consistently ordered $p$-cyclic matrix, for some integer $p \geq 2$; see \cite[ch.~4]{varga1969matrix} for a definition and analysis.  In this paper we are interested in the case $p=2$. 
\begin{definition}
An $n \times n$ matrix $A$ is said to be a consistently ordered $2$-cyclic matrix if it can be symmetrically permuted into the block form
$$
A =
\begin{bmatrix}
D_1 & A_{12}\\[6pt]
A_{21}                    & D_2
\end{bmatrix},
$$
where the matrices $D_1,D_2$ are diagonal and  nonsingular.  
\end{definition}

\noindent Let us define the {\em spectral radius} of a matrix $B$ as
$$ \rho(B) = \max_j |\lambda_j(B)|,$$
where $\sigma(B) := \{ \lambda_j(B)\}_{j=1}^n$ are the eigenvalues of $B$, referred to also as the spectrum of $B$. Then, the asymptotic convergence rate of a stationary iterative method with iteration matrix $B$ is given by
$$R_{\infty}(B) = -\log(\rho(B)).$$

\noindent To state the results related to SOR, let us recall the  {\em Jacobi iteration matrix}:
\begin{equation}
\mathcal J = D^{-1}(E+F) = I - D^{-1} A.
\label{eq:Jac}
\end{equation}
If $A$ is consistently ordered and $2$-cyclic, its Jacobi spectrum is symmetric about the origin \cite[Lemma 2.1]{Young}, that is, $\mu \in \sigma(\mathcal J)$ if and only if $-\mu \in \sigma(\mathcal J)$.  Therefore, if the eigenvalues of $\mathcal J$ are real, then we must have 
$$\sigma(\mathcal J) \subseteq [-\tilde{\mu},\tilde{\mu}] \subseteq \field{R}, \quad \rho(\mathcal J)=\tilde{\mu}.$$
If $\rho(\mathcal J)<1$, the optimal relaxation parameter is famously given by 
\begin{equation} \label{eqn:omega_opt}
\omega_{\rm opt}=\frac{2}{1+\sqrt{1-\rho(\mathcal J)^2}}.
\end{equation}
\noindent
As $\rho(\mathcal J)\rightarrow 1^-$,  by \cite[Cor.~4.9]{varga1969matrix} we have
\begin{equation} \label{eqn:square_root_rates}
R_{\infty}(\mathcal L_{\omega_{\rm opt}}) \sim 2\sqrt{2} \left[R_{\infty}(\mathcal J)\right]^{\frac{1}{2}},
\end{equation}
where we use the standard notation $f(x) \sim g(x)$ as $x \rightarrow a$ to mean $\lim_{x \rightarrow a} \frac{f(x)}{g(x)}=1$. 

\noindent When $A$ is a discretized second-order elliptic partial differential equation (PDE), instead of \eqref{eqn:square_root_rates} estimates are typically of the form 
\begin{equation} \label{eqn:rates_real}
\rho(\mathcal J) \sim 1-\mathcal{O}(h^{2}), \qquad \rho(\mathcal L_{\omega_{\rm opt}}) \sim 1-\mathcal{O}(h), 
\end{equation}
valid for $\rho(\mathcal J) \rightarrow 1^-$, where $h$ is the typical mesh size; see, for example, \cite{FiniteElement, FourierSOR}. 

In this paper, we are interested in generalizing these results to the case where $\tilde{\mu}$ is complex. 
For $\tilde{\mu} \in \field{C}$, we denote by 
$[-\tilde{\mu},\tilde{\mu}]$ the line segment in $\field{C}$ passing through the origin and joining $\pm\tilde{\mu}$, that is,
\begin{equation} \label{eqn:symmetric_line_segment}
[-\tilde{\mu},\tilde{\mu}] =\{ t\tilde{\mu} : t \in [-1,1]\} \subseteq \field{C}, \qquad \operatorname{Re}(\tilde{\mu}) \geq 0,
\end{equation}
where the assumption $\operatorname{Re}(\tilde{\mu}) \geq 0$ is true without loss of generality since at least one of $\pm \tilde{\mu}$ must have a nonnegative real part. 

If $\sigma(\mathcal J) \subseteq [-\tilde{\mu},\tilde{\mu}] \subseteq \field{C}$ and $\pm \tilde{\mu} \in \sigma(\mathcal J)$, it has been shown in \cite{ComplexSOR} that
\begin{equation} \label{eqn:omega_complex_opt}
\omega_{\rm opt}=\frac{2}{1+\sqrt{1-\tilde{\mu}^2}} \in \field{C}, \qquad \rho(\mathcal L_{\omega_{\rm opt}})=|1-\omega_{\rm opt}|.
\end{equation}
However, convergence results of the form \eqref{eqn:square_root_rates} or \eqref{eqn:rates_real} are not provided in \cite{ComplexSOR}.    
There are other derivations of the optimal complex SOR relaxation parameters in the literature for other scenarios; see, for example, \cite{SOR_Ellipse}. However, they do not offer an analysis similar to that presented in this paper.

Our focus on the setting characterized by \eqref{eqn:symmetric_line_segment} is largely motivated by the importance of relevant applications that give rise to this spectral structure.  In particular, this structure often arises in the context of (but is not limited to) linear systems with a complex shift.  There are several examples of this important class of linear systems; see, for example, \cite{hocking2021optimal} and the references therein.  For instance, the damped Helmholtz equation, which we discuss in Section~\ref{sec:examples}, leads to a linear system with a complex shift of the diagonal of a Laplacian \cite{ShiftedLaplacian}.

The remainder of the paper is structured as follows. In Section \ref{sec:conv} we offer a convergence analysis for the case under current consideration, $\tilde{\mu} \in {\field C}$. In Section \ref{sec:examples} we apply our analysis to the damped Helmholtz equation. Finally, in Section~\ref{sec:conc} we draw some conclusions.

\section{Convergence Analysis}
\label{sec:conv}

Our main results follow from the lemma below, which we  state as a stand-alone result in complex analysis.  When we apply it later to SOR, $z \in \field{C}$ will become $\tilde{\mu}$, and $|f(z)|$ will become $\rho(\mathcal L_{\omega_{opt}})$.  We use throughout the principal branch of $\operatorname{Arg}(z)$, that is, $\operatorname{Arg}(z) \in (-\pi,\pi]$.  Similarly, we use the principal branch of the square root, that is, the one with the branch cut on the negative real axis.

\begin{figure}[h!]
\centering
\centering
\begin{tabular}{cc}
\subfloat[$R \leq 1$.]{\includegraphics[width=.48\linewidth]{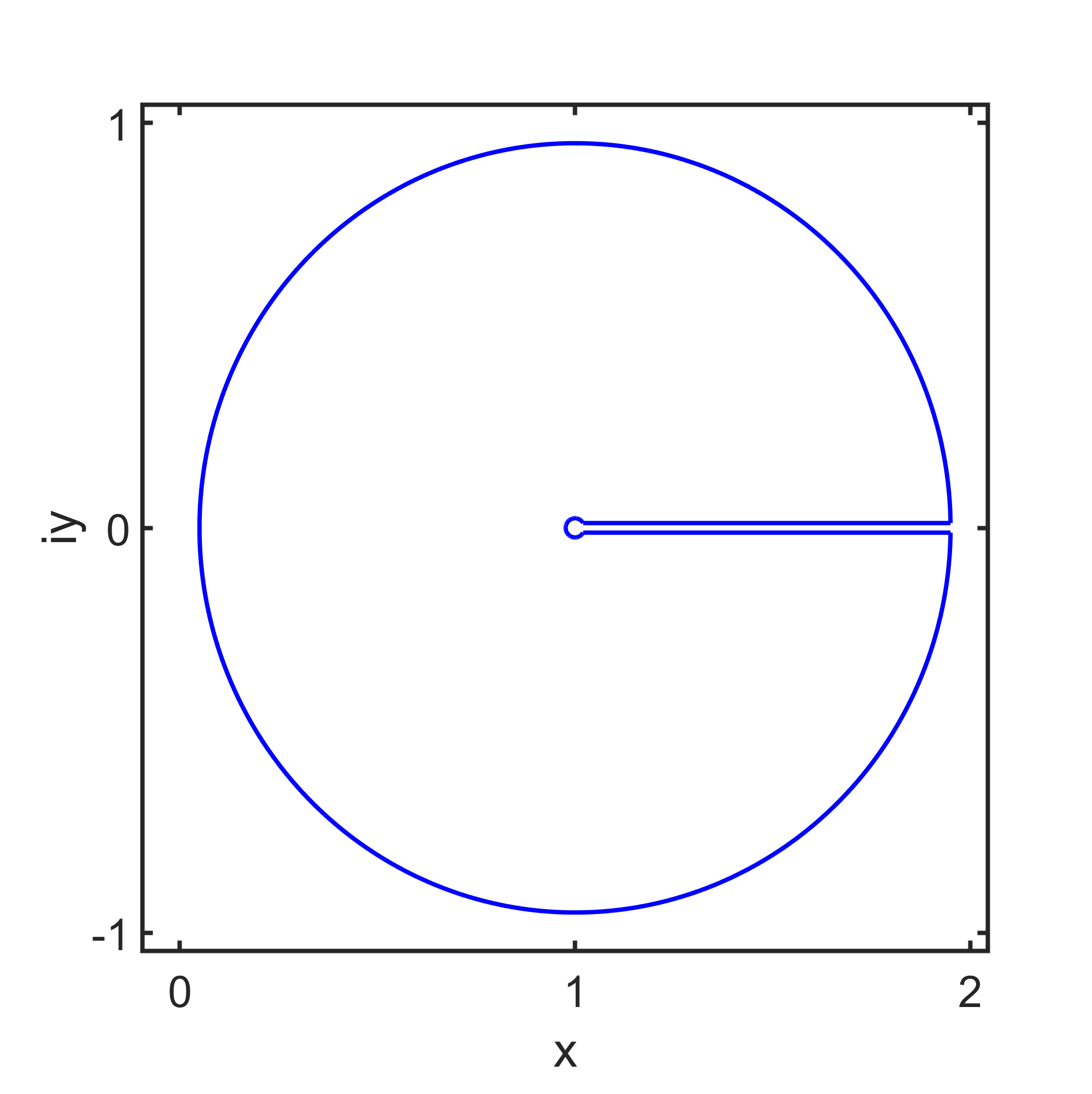}} &
\subfloat[$R>1$.]{\includegraphics[width=.48\linewidth]{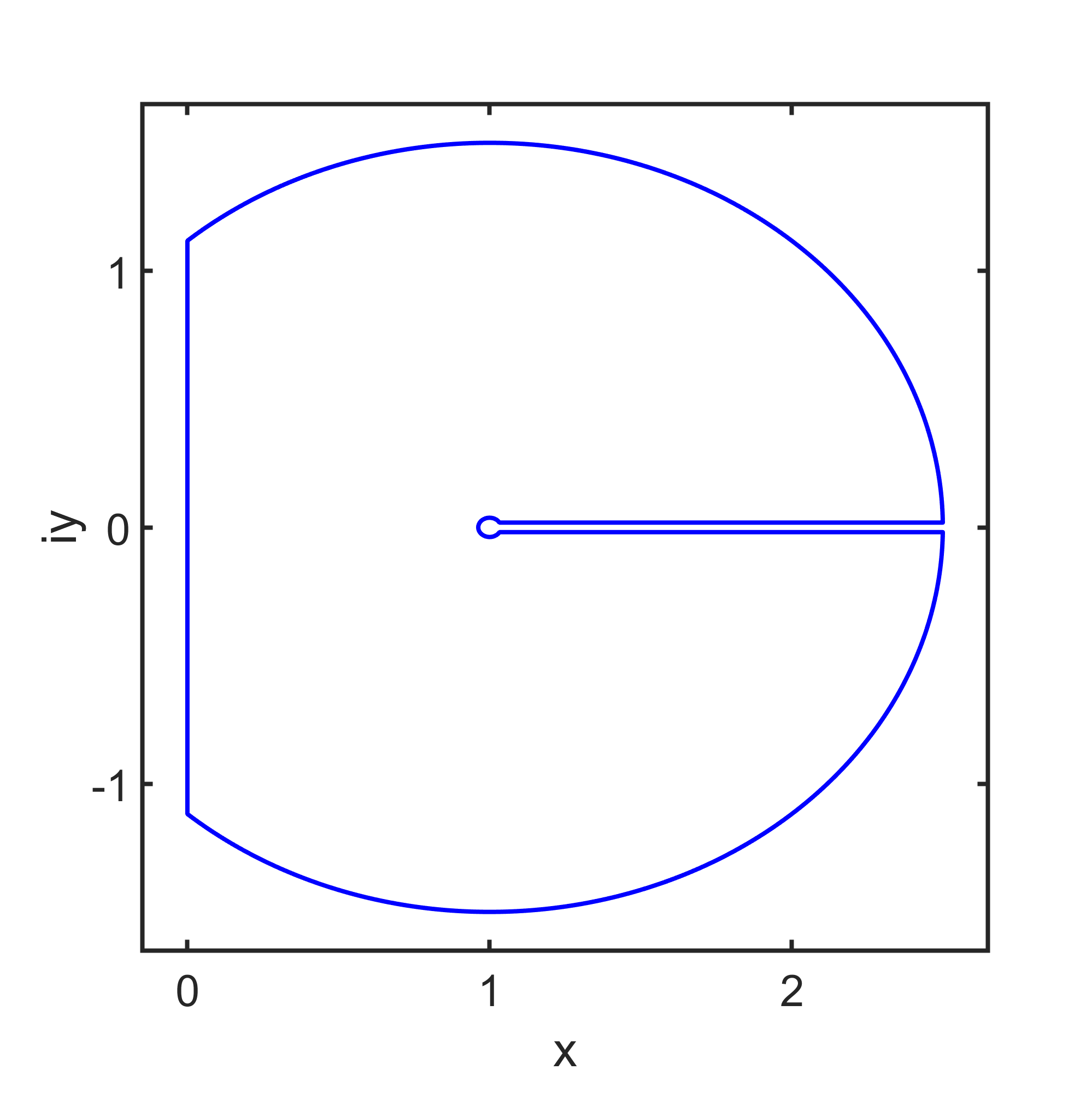}} 
\end{tabular}
\caption[protected]{Illustration of the set $D_R \cap \field{C}_{\geq 0} \setminus [1,\infty)$ from Lemma \ref{lem:f-g-estimate}, for different values of $R$.  For $R>1$, part of the imaginary axis is included in the boundary.}
\label{fig:boundary}
\end{figure}

\begin{lemma}\label{lem:f-g-estimate}
Let
\[
   f(z)=\frac{2}{1+\sqrt{\,1-z^{2}\,}}-1,
   \qquad
   g(z)=1-2\sqrt{2\,(1-z)}.
\]
Define the right half plane $$\field{C}_{\geq 0} =\{ z \in \field{C} : \operatorname{Re}(z) \geq 0\},$$ and the punctured disk $$   D_R=\Bigl\{\,z\in\mathbb C:\;|z-1|<R,\;z\ne1\Bigr\},$$ and let $$\delta = |\operatorname{Arg}(z-1)|.$$ 
Then, for $z\in \field{C}_{\geq 0} \cap D_R$:
\begin{itemize}
    \item[(i)] $|f(z)-1|$ is bounded from above and below as follows:
\begin{equation}\label{eq:main-ineq}
      c_R\,|1-g(z)|
      \;\le\;
      |f(z)-1|
      \;\le\;
      |1-g(z)|,
\end{equation}
where  $0 < c_R$ is a monotonically decreasing function of $R \geq 0$ with $c_0=1$, given explicitly as follows:
\begin{eqnarray}
\label{eq:CR}
c_R & = & \begin{cases}
\displaystyle \frac{1}{\sqrt{2} \left( \frac{1}{\sqrt{2 - R}} + \sqrt{R} \right)} & \text{if } 0 \le R \le 1; \\
\displaystyle \frac{1}{\sqrt{2} \left( \frac{1}{\sqrt{R}} + \sqrt{R} \right)} & \text{if } R > 1.
\end{cases}
\end{eqnarray}
\item [(ii)] $1-|f(z)|$ is bounded from above and below as follows:
 \begin{equation}\label{eq:main-ineq2}
      c^*_R\sin\left(\frac{\delta}{2}+\beta_m\right)|1-g(z)|
      \;\le\;
     1-|f(z)|
      \;\le\;
      \sin\left(\frac{\delta}{2}+\beta_M\right)|1-g(z)|,
\end{equation}
where $\beta_m,\beta_M$ are the values of $0 \leq \beta \leq \frac{1}{2}\arctan|\operatorname{Im}(z)|$ that minimize and maximize $\sin\left(\frac{\delta}{2}+\beta\right)$ respectively, and where
 $0 < c^*_R \leq c_R$ is a monotonically decreasing function of $R \geq 0$ satisfying $c^*_0=c_0=1$,  given explicitly by
 \begin{equation} 
     c^*_R=\frac{c_R}{1+\sqrt{R(2+R)}}.
     \label{eq:CRstar}
     \end{equation}
\end{itemize}

\end{lemma}
\begin{proof}
We proceed by proving claims (i) and (ii) in sequence.
\vskip 2mm
\noindent {\em Proof of claim (i).}
Define  
\[
     p(z):=(1+z)^{-1/2}+(1-z)^{1/2},
\]
and note that $p(z)$ is analytic and non-vanishing on $D_R \cap \field{C}_{\geq 0} \setminus [1,\infty)$.  A few lines of algebra yield
\begin{equation}\label{eq:hF}
     \frac{|f(z)-1|}{|1-g(z)|}
           =\frac1{\sqrt2\,|p(z)|}.
\end{equation}
It follows from the maximum modulus principle that the expression in \eqref{eq:hF} obtains its maximum and minimum values on the boundary of $D_R \cap \field{C}_{\geq 0} \setminus [1,R+1)$, which we illustrate in Figure \ref{fig:boundary}.  
The boundary consists of up to three pieces: the slit $[1,R+1)$, a portion of the imaginary axis (applicable when $R>1$), and all or part of the the circle $|z-1|=R$.  We now consider these three pieces in turn.

We first consider the slit $[1,R+1)$. Let us write $z=x+iy$, where $x, y \in {\mathbb R}$.  One readily computes that $|p(x)|$ is monotonically increasing for $x \in [1,R+1)$, so that \eqref{eq:hF} is maximized at $x=1$, where \eqref{eq:hF} takes on the value $1$.  

Next we consider the imaginary axis in the case $R>1$, where a similar situation occurs. We have
$$|p(iy)|=(1+y^2)^{\frac{1}{4}}+(1+y^2)^{-\frac{1}{4}},$$
which is a monotonically increasing function of $|y|$, hence \eqref{eq:hF} is maximized at $y=0$, where we obtain the value $\frac{1}{2\sqrt{2}}$.  This is smaller than $1$, so it does not contribute to the overall maximum.  At the same time, the minimum values on both the slit and the imaginary axis occur where they intersect the circle $|z-1|=R$.  These minima are therefore absorbed into the next case. 

We now turn to the third and final case of the circle $|z-1|=R$.  Define for convenience
\begin{equation} \label{eqn:s}
    s(z)=\sqrt{1-z^2}.
\end{equation}
On the boundary circle, $|p(z)|$ further simplifies to 
$$|p(z)|=\sqrt{R}\left|1+\frac{1}{s(z)}\right| \quad \mbox{ for } |z-1|=R.$$
Write the boundary point as
$$z=1+Re^{i\theta}$$ and define
$$
      q(\theta)=|2+Re^{i\theta}|=\sqrt{4+R^2+4R\cos\theta}.
$$
Since $\operatorname{Re}(z)\geq 0$, we have $\cos\theta \in [-R^{-1},1]$ and hence $q(\theta) \in[R,2+R].$
Next, observe \(|s(\theta)|=\sqrt{Rq(\theta)}.\)  Since we are working with the principal branch of the square root, we have $\operatorname{Re}(s(\theta)) \geq 0$, from which it follows that
\[
      \sqrt{1+\frac{1}{|s(\theta)|^{2}}} \leq \left|1+\frac{1}{s(\theta)}\right|
      \leq 1+\frac{1}{|s(\theta)|}.
\]
Together, this yields the bound
\[
      \sqrt{R+\frac{1}{q(\theta)}}
      \;\le\;
      \bigl|p(1+Re^{i\theta})|
      \;\le\;
      \sqrt{R}+\frac{1}{\sqrt{q(\theta)}}.
\]
Substituting this in our bounds on $q(\theta)$ gives:
\begin{equation} \label{eqn:not_tight}
      \frac{1}{\sqrt{2}} \leq \sqrt{R+\frac{1}{2+R}}
      \;\le\;
      \bigl|p(1+Re^{i\theta})\bigr|
      \;\le\;
      \sqrt{R}+\frac{1}{\sqrt{R}}. 
\end{equation}
The leftmost inequality gives us an upper bound for \eqref{eq:hF} of $1$, which is the same as the upper bound on the slit and the imaginary axis.  Hence, our overall upper bound is $1$, independent of $R$.  This proves the right-side inequality of \eqref{eq:main-ineq} in claim (i).  

To prove the left-side inequality of \eqref{eq:main-ineq}, note that for $R \leq 1$, a tighter upper bound on $|p(1+Re^{i\theta})|$ than the one in \eqref{eqn:not_tight} may be obtained by observing
$$|p(1+Re^{i\theta})|=\left|\frac{1}{\sqrt{2+Re^{i\theta}}}+\sqrt{-Re^{i\theta}} \right| \leq \frac{1}{\sqrt{2-R}}+\sqrt{R}.$$
Substitution of the above bound for $R \leq 1$ together with the bound from \eqref{eqn:not_tight} for $R>1$ into \eqref{eq:hF} yields the left-side inequality of \eqref{eq:main-ineq} in claim (i) with the explicit expression \eqref{eq:CR}.

\vskip 2mm
\noindent {\em Proof of claim (ii).}
For notational convenience, let us denote henceforth $s(z)$ defined in \eqref{eqn:s}
by $s$. We have
$$f(z)=\frac{1-s}{1+s},\qquad
      |f(z)-1|=\frac{2|s|}{|1+s|}.$$
Consequently,
\begin{equation}\label{eq:num}
   1-|f(z)|=\frac{|1+s|-|1-s|}{|1+s|}
           =\frac{4\,\Rept(s)}{|1+s|\bigl(|1-s|+|1+s|\bigr)},
\end{equation}
where the second equality in \eqref{eq:num} follows after multiplying the top and bottom by $|1+s|+|1-s|$ and from the identity $|1+s|^2-|1-s|^2 = 4\Rept(s)$.  We then have
\begin{equation}\label{eq:ratio'}
      \frac{1-|f(z)|}{|f(z)-1|}
      =\frac{2\,\Rept(s)}{|s|\bigl(|1-s|+|1+s|\bigr)}=\frac{2 \cos(|\operatorname{Arg(s)}|)}{\bigl(|1-s|+|1+s|\bigr)}.
\end{equation}
To find a bound on $|\operatorname{Arg}(s)|$, first note that
$$\operatorname{Arg}(s)=\frac{\pi}{2}+\frac{1}{2}\operatorname{Arg}(z-1)+\frac{1}{2}\operatorname{Arg}(z+1).$$
Next, note that $\operatorname{Re}(z+1)>\operatorname{Re}(z-1)$ together with $\operatorname{Im}(z+1)=\operatorname{Im}(z-1)$ implies that $$0 \leq |\operatorname{Arg}(z+1)|<|\operatorname{Arg}(z-1)| \leq \pi$$ and that $\operatorname{Arg}(z+1)$ has the same sign as $\operatorname{Arg}(z-1)$.  First assume $0 \leq \operatorname{Arg}(z-1) \leq \pi$.  
We obtain
$$|\operatorname{Arg}(z-1)|=\frac{\pi}{2}+\frac{\delta}{2}+\beta(z),$$
where $\delta=|\operatorname{Arg}(z-1)|$ and where $\beta(z)=\frac{1}{2}|\operatorname{Arg}(z+1)|$ obeys the bound
$$0 \leq \beta(z) \leq \frac{1}{2}\arctan\left|\frac{\operatorname{Im}(z)}{1+\operatorname{Re}(z)}\right| \leq \frac{1}{2}\arctan\left|\operatorname{Im}(z)\right|.$$
A similar argument shows that we obtain the same result when $-\pi \leq \operatorname{Arg}(z-1) \leq 0$.
Consequently, we have 
$$\cos(|\operatorname{Arg(s)}|) = \cos\left(\frac{\pi+\delta}{2}+\beta\right)=\sin\left(\frac{\delta}{2}+\beta\right).$$
We define $\beta_m$ and $\beta_M$ to be the values of $0 \leq \beta \leq\frac{1}{2}\arctan\left|\operatorname{Im}(z)\right|$ minimizing and maximizing the above expression, respectively.  
At the same time, we have
$$|s|^2 \leq R(2+R),$$
so that
$$|1-s|+|1+s| \leq 2(1+|s|) \leq 2(1+\sqrt{R(2+R)}).$$
We also have
$$|1-s|+|1+s| \geq |(1-s)+(1+s)|=2.$$
Plugging these bounds into \eqref{eq:ratio'} gives 
\begin{equation}\label{eq:two-sided}
\frac{\sin\left(\frac{\delta}{2}+\beta\right)}{1+\sqrt{R(2+R)}} \leq \frac{1-|f(z)|}{|f(z)-1|} \leq \sin\left(\frac{\delta}{2}+\beta\right).
\end{equation}
Combining \eqref{eq:main-ineq} with \eqref{eq:two-sided} yields \eqref{eq:main-ineq2} of claim (ii), where $c^*_R$ is given in \eqref{eq:CRstar}. 
\end{proof}

\begin{figure}[h!]
\centering
\centering
\begin{tabular}{cc}
\subfloat[Plot of $\frac{1-|f(z)|}{|1-g(z)|}$.]{\includegraphics[width=.48\linewidth]{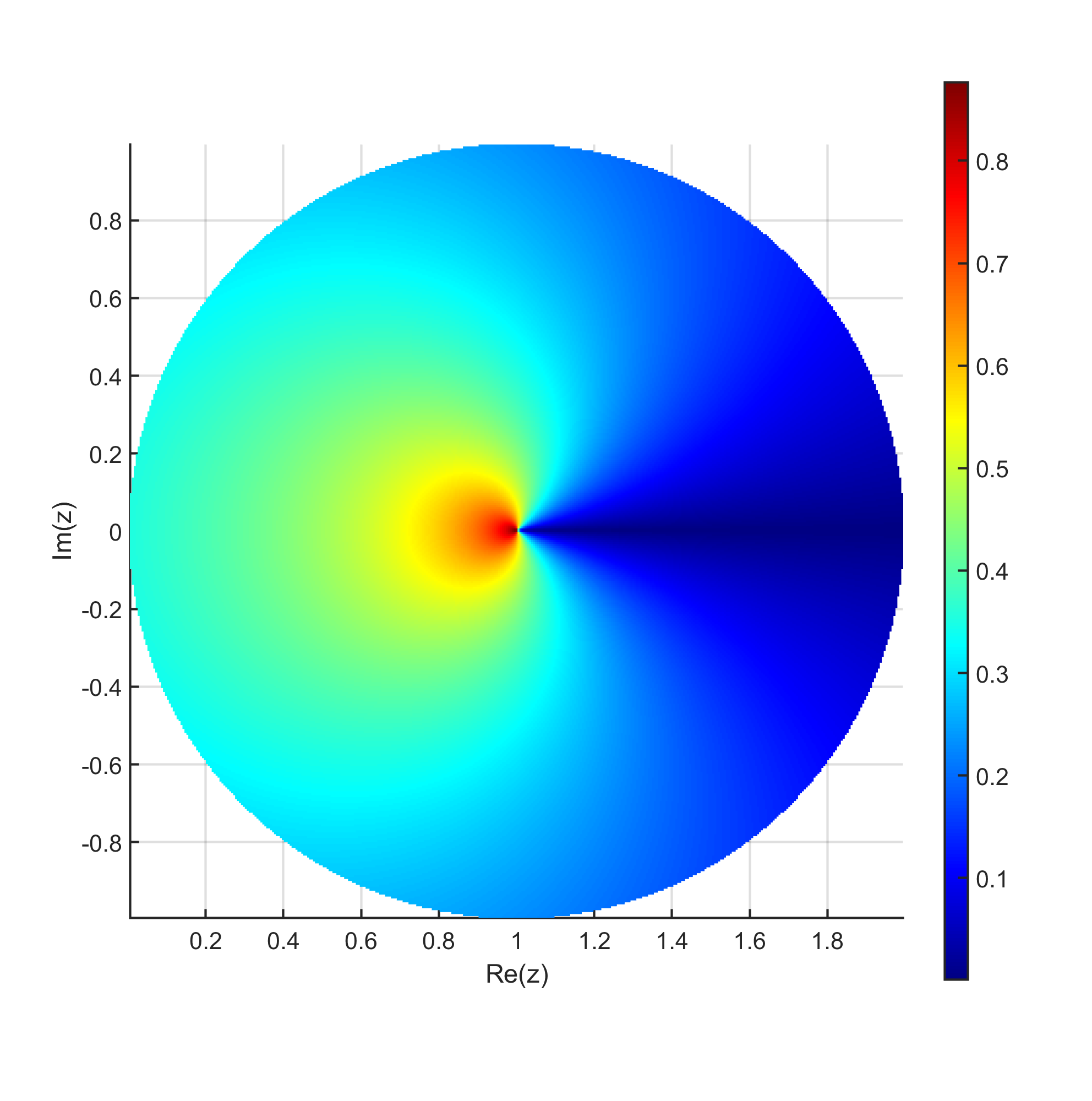}} &
\subfloat[Plot of $|f(z)|$.]{\includegraphics[width=.48\linewidth]{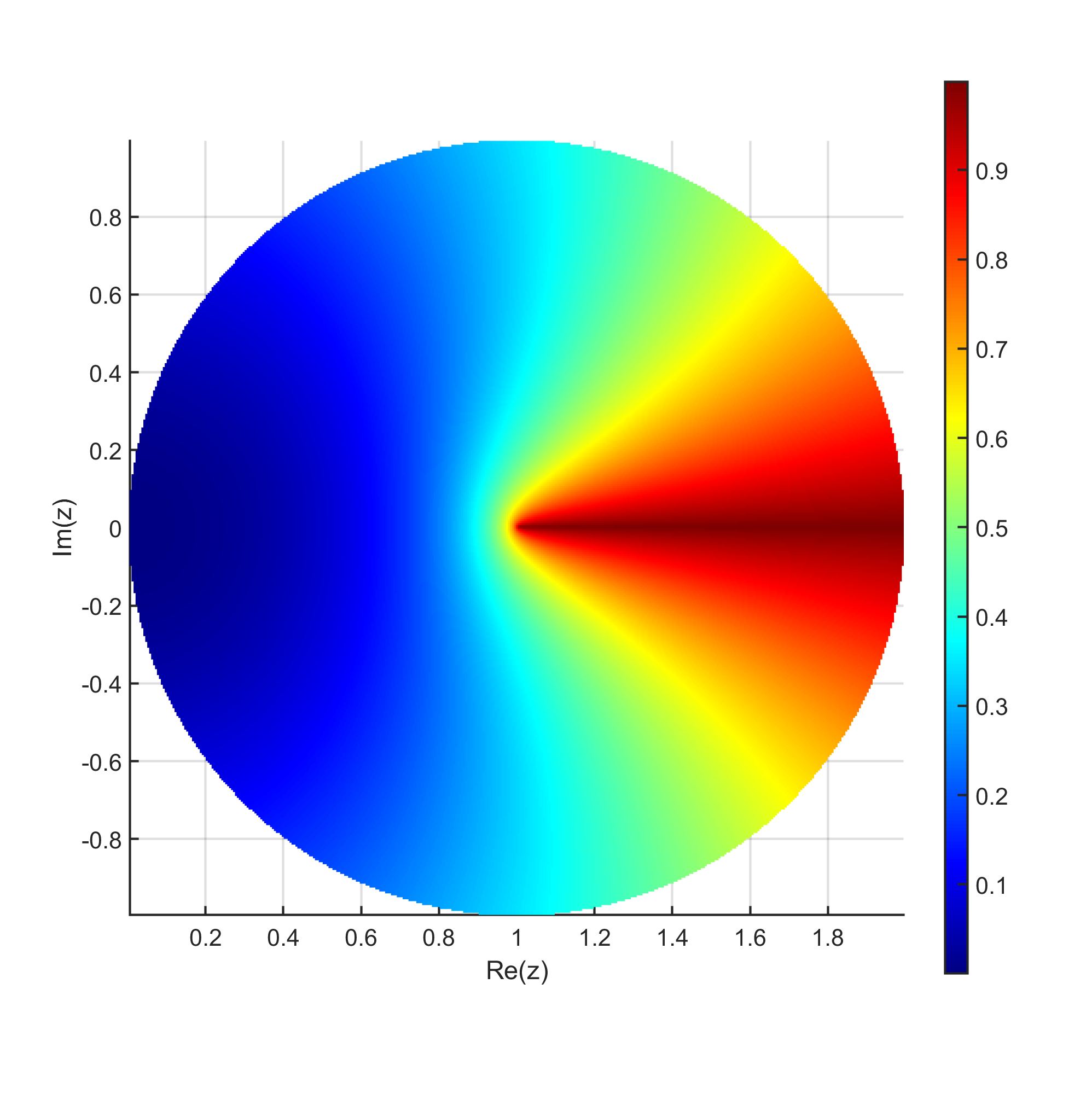}} \\
\end{tabular}
\caption[protected]{Central functions in Lemma \ref{lem:f-g-estimate} and Theorem \ref{thm:main_rates}, for $R \le 1$.}
\label{fig:bound}
\end{figure}

\noindent Fig. \ref{fig:bound}(a) provides a plot of $\frac{1-|f(z)|}{|1-g(z)|}$ for $f$ and $g$ as defined in Lemma \ref{lem:f-g-estimate}, with $R \leq 1$.  We observe a strong sensitivity to $|\operatorname{Arg}(z-1)|$, especially close to $z=1$, and the ratio vanishes on the slit $[1,\infty)$.  This is reflected in the bounds provided by Lemma \ref{lem:f-g-estimate}, which become tight as $z \rightarrow 1$ and on the slit $[1,\infty)$.  Figure \ref{fig:bound}(b) provides a similar plot of $|f(z)|$ where we observe $|f(z)| \leq 1$ with equality if and only if $z \in [1,\infty)$, a statement we will prove in Theorem \ref{thm:main_rates}.

\begin{remark}
Lemma \ref{lem:f-g-estimate} says that near the point $z=1$, we have $$|f(z)-1| \approx 2\sqrt{2}\sqrt{ |1-z|}.$$
The fast convergence of SOR is due to the square root above, which is a direct consequence of the breakdown of analyticity at the point $z=1$.
Indeed, if $f(z)$ were analytic at $z=1$, then we would have, for $z \approx 1$, 
$|f(z)-1| \approx |f'(1)||z-1|$ by the Taylor expansion, which would in turn break the SOR ``magic.''  
\end{remark}

\noindent We are now ready to prove the main result of this section, and indeed this paper.
\begin{theorem} \label{thm:main_rates}
Given a linear system  \eqref{eq:system} with a matrix $A$ that is consistently ordered and 2-cyclic, consider the  splitting \eqref{eq:splitting} and assume that the Jacobi spectrum associated with the iteration matrix \eqref{eq:Jac} satisfies $\rho(\mathcal J) \subseteq [-\tilde{\mu},\tilde{\mu}] \subseteq \field{C}$, as well as $\pm \tilde{\mu} \in \sigma(\mathcal J)$. Using the convention $\operatorname{Re}(\tilde{\mu}) \geq 0$ from \eqref{eqn:symmetric_line_segment},  define 
$$\delta = |\operatorname{Arg}(\tilde{\mu}-1)|.$$
Then, if \eqref{eq:system} is solved using SOR as per \eqref{eq:SOR}, the optimal relaxation parameter $\omega_{\rm opt}$ is given by the formula on the left of \eqref{eqn:omega_complex_opt},
and the SOR spectral radius $\rho(\mathcal L_{\omega_{\rm opt}})$ on the right of \eqref{eqn:omega_complex_opt} obeys
\begin{equation} \label{eqn:double_sided_bound}
2\sqrt{2}c^*_R\sin\left(\frac{\delta}{2}+\beta_m\right)\sqrt{|1-\tilde{\mu}|} \leq 1-\rho(\mathcal L_{\omega_{\rm opt}}) \leq 2\sqrt{2} \sin\left(\frac{\delta}{2}+\beta_M\right)\sqrt{|1-\tilde{\mu}|} ,  
\end{equation}
where $\beta_m,\beta_M$ are the values of $0 \leq \beta \leq \frac{1}{2}\arctan|\operatorname{Im}(\tilde{\mu})|$ that minimize and maximize $\sin\left(\frac{\delta}{2}+\beta\right)$ respectively, and where $0 < c^*_R$ is the monotonically decreasing function of $R \geq 0$ obeying $c^*_0=1$ and given explicitly by \eqref{eq:CRstar}. 
Consequently, we have 
$$\rho(\mathcal L_{\omega_{\rm opt}}) < 1 \quad {\rm if\ }\tilde{\mu} \in \field{C}_{\geq 0} \setminus [1,\infty)$$
and $$\rho(\mathcal L_{\omega_{\rm opt}}) = 1 \quad {\rm if\ }\tilde{\mu} \in [1,\infty).$$
\end{theorem}
\begin{proof}
First, note that $|1-g(z)|=2\sqrt{2}\sqrt{|1-z|}$.  Taking $z=\tilde{\mu}$, inequality \eqref{eqn:double_sided_bound} 
follows immediately from \eqref{eq:main-ineq2} in Lemma \ref{lem:f-g-estimate}.  This inequality becomes an equality if and only if $\tilde{\mu} \in [1,\infty)$, as $\sin\left(\frac{\delta}{2}+\beta_m\right)=\sin\left(\frac{\delta}{2}+\beta_M\right)=0$ if and only if $\tilde{\mu} \in [1,\infty)$.  It follows that $\rho(\mathcal L_{\omega_{\rm opt}}) \leq 1$ with equality if and only if $\tilde{\mu}$ is on the slit $[1,\infty)$.

\end{proof}

\noindent It is illuminating to consider the case $\tilde{\mu} \approx 1$, or equivalently $R \ll 1$, as follows.
\begin{corollary} \label{col:generalize}
Given the assumptions of Theorem \ref{thm:main_rates}, we have
\begin{equation}
    \label{eq:R}
    R_{\infty}(\mathcal L_{\omega_{\rm opt}}) \sim 2\sqrt{2} \sin\left(\frac{\delta}{2}\right)\sqrt{|1-\tilde{\mu}|}, \qquad \tilde{\mu} \rightarrow 1, \qquad \tilde{\mu} \in \field{C}_{\geq 0} \setminus [1,\infty).
\end{equation}
\end{corollary}
\begin{proof}
Assume $R \leq 1$.  By \eqref{eqn:double_sided_bound} we have for all $\tilde{\mu} \in D_R \setminus [1,\infty)$
$$
\frac{c^*_R\sin\left(\frac{\delta}{2}+\beta_m\right)}{\sin\left(\frac{\delta}{2}\right)} \leq \frac{1-\rho(\mathcal L_{\omega_{\rm opt}})}{2\sqrt{2}\sin\left(\frac{\delta}{2}\right)\sqrt{|1-\tilde{\mu}|}} \leq \frac{\sin\left(\frac{\delta}{2}+\beta_M\right)}{\sin\left(\frac{\delta}{2}\right)}.$$
As $R \rightarrow 0$ we have $c^*_R \rightarrow 1$ while $\beta_m,\beta_M \rightarrow 0$, so that both the upper and lower bounds in the inequality tend to $1$, yielding
$$1-\rho(\mathcal L_{\omega_{\rm opt}}) \sim 2\sqrt{2}\sin\left(\frac{\delta}{2}\right)\sqrt{|1-\tilde{\mu}|}, \qquad \tilde{\mu} \rightarrow 1, \qquad \tilde{\mu} \in \field{C}_{\geq 0} \setminus [1,\infty).$$
At the same time, it follows from l'H\^opital's rule that
$$1-\rho(\mathcal L_{\omega_{\rm opt}}) \sim -\log(\rho(\mathcal L_{\omega_{\rm opt}})), \qquad \rho(\mathcal L_{\omega_{\rm opt}}) \rightarrow 1^-.$$
The desired result then follows from the observation that $\rho(\mathcal L_{\omega_{\rm opt}}) \rightarrow 1^-$ as $\tilde{\mu} \rightarrow 1$ with $\tilde{\mu} \in \field{C}_{\geq 0} \setminus [1,\infty)$ and the product rule for limits.
\end{proof}
\begin{remark}
    The result stated in Corollary \ref{col:generalize}, namely Eq. \eqref{eq:R}, reduces to the classical result \eqref{eqn:square_root_rates} of Varga \cite[p. 126]{varga1969matrix} when $\sigma(J) \subseteq (-1,1)$, as in that case we have $\sin\left(\frac{\delta}{2}\right)=\sin\left(\frac{\pi}{2}\right)=1$ and $\tilde{\mu}=\rho(J)<1$, so that $\sqrt{|1-\tilde{\mu}|}=[R_{\infty}(\mathcal J)]^{\frac{1}{2}}$.  Moreover, our restriction of avoiding the slit $[1,\infty)$ when taking the limit in the complex plane reduces to Varga's requirement that we must have $\rho(J)\rightarrow 1^-$.
\end{remark}

\section{Example: damped Helmholtz equation in two dimensions} \label{sec:examples}

Consider a two-dimensional model damped Helmholtz equation 
\begin{equation} \label{eqn:shifted_laplacian}
-\Delta u - (1-i\alpha)k^2(\vec{x})u = f(\vec{x}),
\end{equation}
where $\alpha>0$ is a damping parameter.  Often, a discrete version of this equation arises also in the context of preconditioning of the numerical solution of {\em undamped} Helmholtz equation \cite{ShiftedLaplacian}.

We discretize \eqref{eqn:shifted_laplacian} on the the unit square, $(0,1) \times (0,1)$, with Dirichlet boundary conditions. For simplicity, we assume that the wavenumber $k(\vec{x})=k$ is constant. We use a uniform mesh with $N$ meshpoints in each direction, and the standard 5-point stencil is used for the Laplacian.  The meshsize is given by $h = \frac{1}{N+1},$
and the dimensions of the corresponding matrix are $N^2 \times N^2$.  

We run our experiments in a {\textsc {Matlab}} environment on a standard laptop. We test with matrices of dimensions varying from $6,400 \times 6,400$ ($N=80$) to $409,600 \times 409,600$ ($N=640$). We take a zero  initial guess and a random  right-hand side vector, and we stop the iteration once the relative residual norm has reached $10^{-6}$. We run our experiments for $\alpha=\frac12, \frac14, \frac18, \frac{1}{16}$, and repeat them for two values of $k$ chosen so that $kh \leq \frac{\pi}{5}$ on all grids, to avoid pollution effects; see \cite{bayliss1985accuracy, hocking2021optimal}. 
Our results are shown in Tables \ref{tab:iterations} and \ref{tab:iterations2}.

The Jacobi spectrum $\sigma(\mathcal J)$ is given by the eigenvalues
$$ \lambda_{j,\ell}=  \frac{2 (\cos\left (j\pi h\right) + \cos\left (\ell \pi h\right ))}{4-(1-i\alpha)k^2h^2}, \qquad 1 \leq j, \ell \leq N,$$
and it is contained in the symmetric line segment $[-\tilde{\mu},\tilde{\mu}]$, where
$$\tilde{\mu}=\frac{\cos(\pi h)}{1-\gamma h^2}, \qquad \gamma = \frac{(1-i\alpha)k^2}{4} \in \field{C}.$$
At the same time, $\tilde{\mu} \in \sigma(\mathcal J)$.  Hence, $\omega_{\rm opt}$ is given by \eqref{eqn:omega_complex_opt}, and the analysis of Section \ref{sec:conv} is applicable.
Expanding the numerator and denominator of $\tilde{\mu}$ in a Taylor series, we obtain
$$\tilde{\mu}=\left(1-\frac{\pi^2h^2}{2}+\mathcal{O}(h^4)\right)\left(1+\gamma h^2+\mathcal{O}(h^4)\right)=1+\left(\gamma-\frac{\pi^2}{2}\right)h^2+\mathcal{O}(h^4), $$
from which it follows that
$$
|\tilde{\mu}-1| = \left(\frac{k^2}{4}\sqrt{\alpha^2 + \left(1-\frac{2\pi^2}{k^2}\right)^2}\right)h^2+\mathcal{O}(h^4) $$
and 
$$|\operatorname{Arg}(\tilde{\mu}-1)| = \arctan\left| \frac{\alpha}{1-\frac{2\pi^2}{k^2}}\right|+\mathcal{O}(h^4).
$$
If $k^2 \gg 2\pi^2$ and $\alpha$ is small, this is approximately
\begin{equation*}
|\tilde{\mu}-1| \approx \frac{k^2}{4}h^2 \quad {\rm and} \quad
|\operatorname{Arg}(\tilde{\mu}-1)| \approx \alpha.
\end{equation*}
Substituting the above expressions into \eqref{eq:R} from Corollary \ref{col:generalize}, we obtain
\begin{equation} \label{eqn:helmholtz_scaling}
R_{\infty}(\mathcal L_{\omega_{\rm opt}}) \sim \frac{\alpha k h}{\sqrt{2}}, \qquad \alpha, h \ll 1, \quad k^2 \gg 2\pi^2.
\end{equation}

\begin{table}[ht]
\centering
\begin{tabular}{|c|cccc|}
\hline
$\alpha$ & $N=80$ & $N=160$ & $N=320$ & $N=640$\\ \hline
$\nicefrac{1}{2}$ & 114 & 230 & 461 & 926\\
$\nicefrac{1}{4}$ & 167 & 330 & 656 & 1305\\
$\nicefrac{1}{8}$ & 275 & 555 & 1123 & 2295\\
$\nicefrac{1}{16}$ & 507 & 1029 & 2065 & 4150\\
\hline
\end{tabular}
\caption{Number of iterations required for convergence of SOR with optimal complex $\omega$ applied to the discretization of \eqref{eqn:shifted_laplacian}, with $k=16\pi$ fixed and for various $\alpha$ and grid sizes $N$.}
\label{tab:iterations}
\end{table}

\begin{table}[ht]
\centering
\begin{tabular}{|c|cccc|}
\hline
$\alpha$ & $N=80$ & $N=160$ & $N=320$ & $N=640$\\ \hline
$\nicefrac{1}{2}$ & 167 & 332 & 659 & 1312\\
$\nicefrac{1}{4}$ & 276 & 564 & 1173 & 2428\\
$\nicefrac{1}{8}$ & 514 & 1032 & 2079 & 4195\\
$\nicefrac{1}{16}$ & 1013 & 2018 & 4082 & 8316\\
\hline
\end{tabular}
\caption{Number of iterations required for convergence of SOR with optimal complex $\omega$ applied to the discretization of \eqref{eqn:shifted_laplacian}, with $k=8\pi$ fixed and for various $\alpha$ and grid sizes $N$.}
\label{tab:iterations2}
\end{table}

In this regime, we therefore expect $\mathcal{O}(\alpha^{-1}k^{-1}N)$ iterations for convergence.  
As such, we expect the entries within a row of Tables \ref{tab:iterations} and \ref{tab:iterations2}  to double as we move from left to right, and the entries within a column to double as we move from top down.  Similarly, we expect the entries of Table \ref{tab:iterations2} to be roughly double those of Table \ref{tab:iterations}.  For small values of $\alpha$, this is what we see; for larger $\alpha$ we notice some modest deviation.  This is to be expected, as the assumption $\alpha \ll 1$ no longer holds.

\noindent Equation \eqref{eqn:helmholtz_scaling}
 predicts constant iterations for constant $kh$.  To test this, we note that $kh$ in the $m$th column of Table \ref{tab:iterations} is the same as in the $(m-1)$st column of Table \ref{tab:iterations2}.  As before, we observe modest deviations for larger $\alpha$, but this becomes asymptotically true as $\alpha \rightarrow 0$.

\begin{figure}[h!]
\centering
\begin{tabular}{cc}
\subfloat[$\alpha = \nicefrac{1}{2}$.]{\includegraphics[width=.46\linewidth]{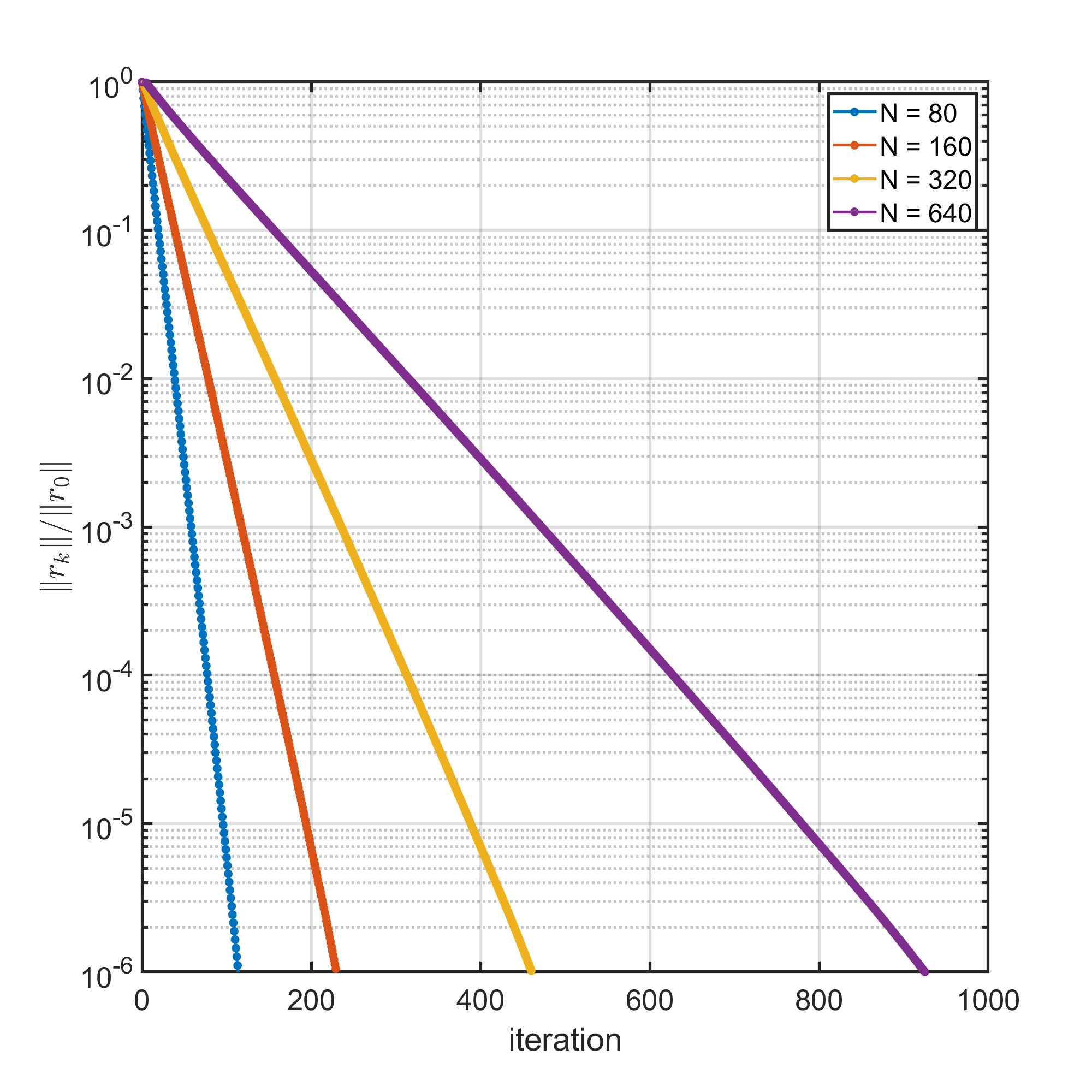}} &
\subfloat[$\alpha = \nicefrac{1}{16}$.]{\includegraphics[width=.46\linewidth]{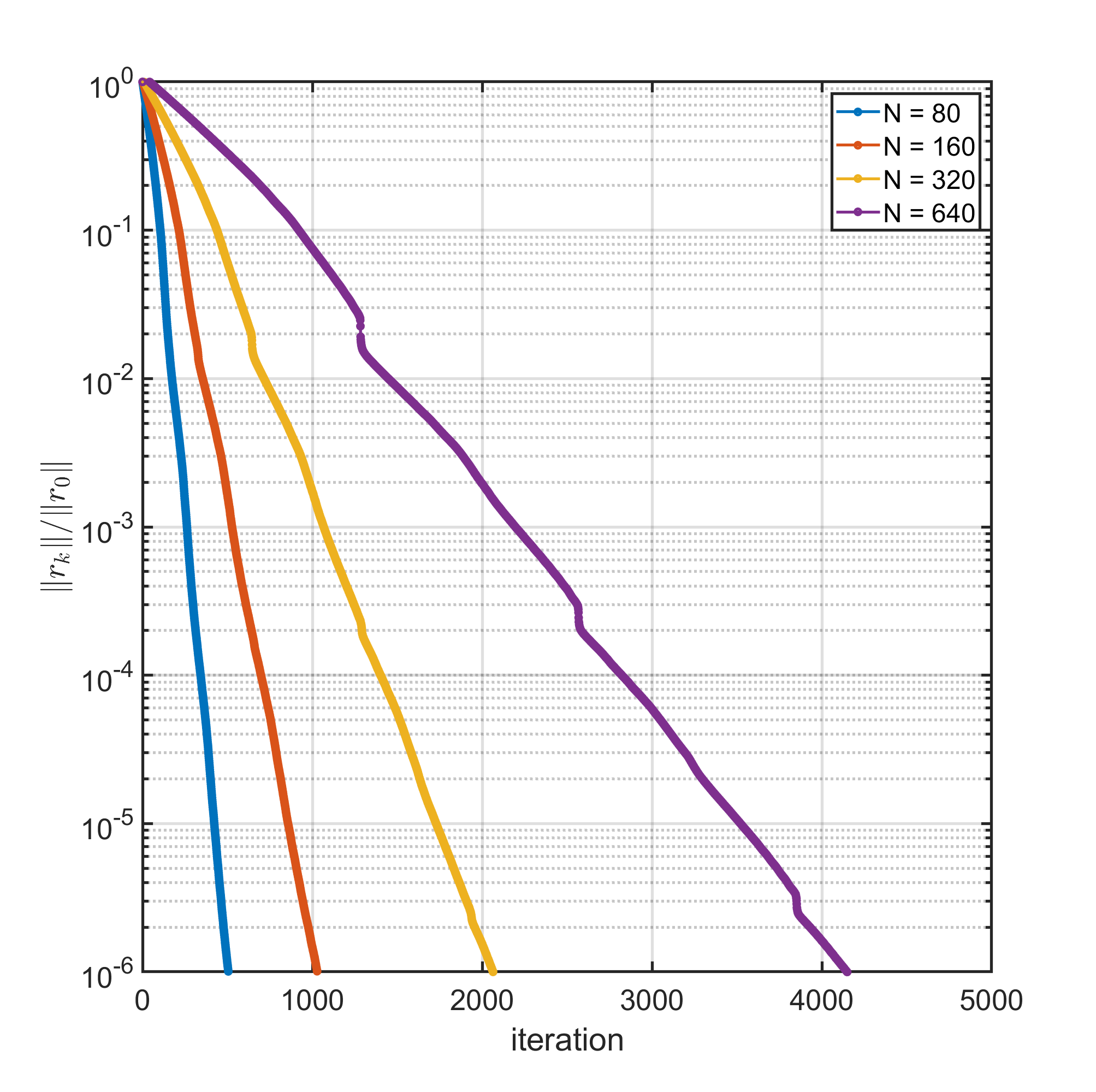}} \\
\end{tabular}
\caption[protected]{Convergence graphs of SOR with optimal complex $\omega$ applied to the discretization of \eqref{eqn:shifted_laplacian}, with $k=16\pi$ and two different values of $\alpha$.  
 }
\label{fig:convergence_curves}
\end{figure}

Figure \ref{fig:convergence_curves} gives some example convergence curves for large and small $\alpha$, with $k=16\pi$ fixed. For large $\alpha$, the convergence curve appears perfectly linear when viewed on a log scale.  For a smaller $\alpha$ there are some departures from linearity, but they are mild, and we do not observe any behavior that one sees for highly non-normal matrices, such as an increase or a slow decrease in the norm of the residual in earlier iterations, before entering the asymptotic regime.

\section{Concluding remarks}
\label{sec:conc}

We have extended classical results of Young and Varga \cite{Young, varga1969matrix} on the relationship between $\rho(\mathcal J)$ and $\rho(\mathcal L_{\omega_{\rm opt}})$ for $\sigma(\mathcal J) \subseteq [-\tilde{\mu},\tilde{\mu}] \subseteq (-1,1)$ as $\tilde{\mu} \rightarrow 1^-$ to the case where $\tilde{\mu} \in \field{C}_{\geq 0}$.  In the real case, one finds $$1-\rho(\mathcal L_{\omega_{\rm opt}}) \sim 2\sqrt{2}[1-\rho(\mathcal J)]^{\frac{1}{2}} = {2\sqrt{2}}\sqrt{|\tilde{\mu}-1|},$$ where the right-hand equality follows since $\tilde{\mu}$ and $\rho(\mathcal J)$ are equal and may be used interchangeably.  Our interest is in the complex case, in which formula \eqref{eqn:omega_complex_opt} was derived for the optimal parameter in \cite{ComplexSOR}. 
Similarly to the real case, we find that $1-\rho(\mathcal L_{\omega_{\rm opt}})$ scales like $\sqrt{|\tilde{\mu}-1|}$.  On the other hand, unlike the real case,  $\tilde{\mu} \in \field{C}$ may now approach $1$ from an infinite number of directions, and indeed we find that $\operatorname{Arg}(\tilde{\mu}-1)$ also has a major impact on convergence, with convergence rates going to zero as $\operatorname{Arg}(\tilde{\mu}-1) \rightarrow 0$.

Since the emergence of Krylov subspace solvers \cite{saad2003iterative} as the gold standard of modern iterative methods, SOR is no longer considered state-of-the-art as a stand-alone solver, but it nonetheless remains an important building block for modern solvers. Notably, it is used as a smoother in the context of multigrid.  Moveover, there is overlap in the theory of SOR in its different roles. For example, in our previous work \cite[Proposition~5.2]{hocking2021optimal}, a direct mathematical connection was established between the Local Fourier Analysis (LFA) smoothing factor of SOR as a smoother and its spectral radius as a solver in some situations. 

A natural direction for future research is the consideration of the Jacobi spectrum other than a line segment, and in what circumstances SOR's speed of convergence can be maintained.

\bibliography{references}
\bibliographystyle {plain}

\end{document}